\documentclass[11pt, a4paper, twoside]{article}
\usepackage{amsfonts}
\usepackage{mathrsfs}
\usepackage{amsmath}
\usepackage{amssymb}
\usepackage{fancyhdr}
\usepackage{hyperref}

\numberwithin{equation}{section}

\allowdisplaybreaks

\begin{document}

\fancyhf{}

\fancyhead[EC]{Heping Liu and Manli Song}

\fancyhead[EL]{\thepage}

\fancyhead[OC]{A functional calculus and restriction theorem on H-type groups}

\fancyhead[OR]{\thepage}

\renewcommand{\headrulewidth}{0pt}
\renewcommand{\thefootnote}{\fnsymbol {footnote}}

\title{\textbf{A functional calculus and restriction theorem on H-type groups}}
\footnotetext {{}\emph{Key words and phrases}: H-type group, scaled special Hermite expansion, joint functional calculus, restriction operator.}
\setcounter{footnote}{0}
\author{Heping Liu and Manli Song
\footnote {Corresponding author.}}
\footnotetext{The first author is supported by National Natural Science Foundation of
China under Grant \#11371036 and the Specialized Research Fund for the Doctoral Program of Higher
Education of China under Grant \#2012000110059. The second author is supported by the China Scholarship Council under Grant \#201206010098.}
\date{}
\maketitle

\begin{abstract}
Let $L$ be the sublaplacian and $T$ the partial Laplacian with respect to central variables on H-type groups.
We investigate a class of invariant differential operators by the joint functional calculus of $L$ and $T$.
We establish Stein-Tomas type restriction theorems for these operators. In particular, the asymptotic behaviors of restriction estimates are given.
\end{abstract}

\newtheorem{theorem}{Theorem}[section]
\newtheorem{preliminaries}{Preliminaries}[section]
\newtheorem{definition}{Difinition}[section]
\newtheorem{main result}{Main Result}[section]
\newtheorem{lemma}{Lemma}[section]
\newtheorem{proposition}{Proposition}[section]
\newtheorem{corollary}{Corollary}[section]
\newtheorem{remark}{Remark}[section]

\section[Introduction]{Introduction}
\quad The restriction theorem for the Fourier transform plays an important role in harmonic analysis as well as in the theory of partial differential equations. The initial work on restriction theorem was given by E. M. Stein \cite{S} that the Fourier transform of an $L^p$-function on $\mathbb{R}^n$ has a well-defined restriction to the unit sphere $S^{n-1}$ which is square integral on $S^{n-1}$. The result is listed as follows:
\begin{theorem}(Stein-Tomas) \indent Let $1 \leq p \leq \frac{2n+2}{n+3}$. Then the estimate
\begin{equation}
||\hat{f}||_{L^2(S^{n-1})} \leq C ||f||_{L^p(\mathbb{R}^n)}         \label{equ:Stein-Tomas}
\end{equation}
\end{theorem}
holds for all functions $f \in L^p(\mathbb{R}^n)$.

A simple duality argument shows that the estimate \eqref{equ:Stein-Tomas} is equivalent to the following estimate:
\begin{equation}
||f*\widehat{d\sigma_r}||_{p^{'}} \leq C_r ||f||_p \label{equ:dual}
\end{equation}
for all Schwartz functions $f$ on $\mathbb{R}^n$, where $\frac{1}{p}+\frac{1}{p'}=1$ and $d\sigma_r$ is the surface measure on the sphere with radius r.

Moreover, according to the Knapp example \cite{S}, the above estimates \eqref{equ:Stein-Tomas} and \eqref{equ:dual} fail if $\frac{2n+2}{n+3}<p\leq2$.

Later, many authors have worked on the topic and various new restriction theorems have been proved. The study of restriction theorems has recently obtained more and more attention. Recent progress on the restriction theorems can be found in \cite{Tao}. To generalize the restriction theorem on Heisenberg group, D. M$\ddot{u}$ller \cite{M} established the boundedness of the restriction operator with respect to the mixed $L^p-$norm and also give a counterexample to show that the estimate between Lebesgue spaces for the restriction operator is necessarily trivial, due to the fact that the center of the Heisenberg group is of dimension 1.

On an H-type group, let $T$ be the Laplacian on the centre and $L$ the sublaplacian. It is well known that $L$ is positive and essentially self-adjoint. Let $L=\int_0 ^{+\infty} \lambda dE(\lambda)$ be the spectral decomposition of $L$. Then the restriction operator can be formally written $\mathcal{P}_\lambda f=\delta_\lambda(L)f=\underset{\epsilon \to +\infty}{\lim}\chi_{(\lambda -\epsilon,\lambda+\epsilon)}(L)f$ which is well defined for a Schwartz function $f$, where $\chi_{(\lambda -\epsilon,\lambda+\epsilon)}$ is the characteristic function of the interval $(\lambda -\epsilon, \lambda+ \epsilon)$. Liu and Wang \cite{LW} investigated the restriction theorem for the sublaplacian $L$ on H-type groups with the center whose dimension is more than 1. They give the following result:
\begin{theorem} \label{Liu and Wang}\indent Let $G$ be an H-type group with the underlying manlifold $\mathbb{R}^{2n+m}$, where $m>1$ is the dimension of the center. Suppose $1 \leq p \leq \frac{2m+2}{m+3}$. Then the following estimate
\begin{equation*}
||\mathcal{P}_\lambda f||_{p'} \leq C \lambda^{2(n+m)(\frac{1}{p}-\frac{1}{2})-1}||f||_p, \lambda>0
\end{equation*}
holds for all Schwartz functions $f$ on $G$.
\end{theorem}

In a recent paper, V. Casarino and P. Ciatti \cite{CC} extend the work of M$\ddot{u}$ller, Liu and Wang on M$\acute{e}$tivier groups. They prove the following result.
\begin{theorem}\label{CC} \indent Let $G$ be a M$\acute{e}$tivier group, with Lie algebra $\mathfrak{g}$. Let $\mathfrak{z}$ and $\mathfrak{v}$ denote, respectively, the centre of $\mathfrak{g}$ and its orthogonal complement.\\
If $dim\mathfrak{z}=d$ and $dim\mathfrak{v}=2n$, and if $1\leq r\leq \frac{2d+2}{d+3}$, then for all $p, q$ satisfying $1\leq p\leq2\leq q\leq\infty$ and for all Schwartz functions $f$, we have
\begin{equation*}
||\mathcal{P}_\lambda f||_{L^{r'}(\mathfrak{z})L^q(\mathfrak{v})}\leq C\lambda^{d(\frac{2}{r}-1)+n(\frac{1}{p}-\frac{1}{q})-1}||f||_{L^r(\mathfrak{z})L^p(\mathfrak{v})}, \lambda>0
\end{equation*}
\end{theorem}

Although V. Casarino and P. Ciatti \cite{CC} investigate the joint functional calculus of $L$ and $T$ on the Heisenberg group, they only prove a restriction theorem for the sublaplacian $L$ on M$\acute{e}$tivier groups. The invariant differential operators related to the joint functional calculus of $L$ and $T$ on H-type groups do not have the homogeneous properties in general. Thus the asymptotic behaviors of restriction estimates for these operators are also interested. In this article we will show restriction theorems for these operators on H-type groups. The range of $p$ depends on the dimension of the center. In particular, the asymptotic behaviors of restriction estimates are given.

The outline of the paper is as follows. In the second section, we provide the necessary background for the H-type group. In the next section, by introducing the joint functional calculus of $L$ and $T$, the restriction operator can be computed explicitly. In the fourth section, we prove the restriction theorem on H-type groups. Finally, in the last section, we show that the range of $p$ in the restriction theorem is sharp.

\section[Preliminaries]{Preliminaries}
\begin{bf} {Definition 2.1}\end{bf}. (H-type Group)\indent Let $\mathfrak{g}$ be a two step nilpotent Lie algebra endowed with an inner product $\langle \cdot,\cdot \rangle$. Its center is denoted by $\mathfrak{z}$. $\mathfrak{g}$ is said to be of H-type if $[\mathfrak{z}^{\bot},\mathfrak{z}^{\bot}]=\mathfrak{z}$ and for every $t \in \mathfrak{z}$, the map $J_t: \mathfrak{z}^{\bot} \rightarrow \mathfrak{z}^{\bot}$ defined by
\begin{equation*}
\langle J_t u, w \rangle:=\langle t, [u,w] \rangle, \forall u, w \in \mathfrak{z}
\end{equation*}
is an orthogonal map whenever $|t|=1$.

An H-type group is a connected and simply connected Lie group $G$ whose Lie algebra is of H-type.

For a given $0 \neq a \in \mathfrak{z}^*$, the dual of $\mathfrak{z}$, we can define a skew-symmetric mapping $B(a)$ on $\mathfrak{z}^{\bot}$ by
\begin{equation*}
\langle B(a)u,w \rangle =a([u,w]), \forall u,w \in \mathfrak{z}^{\bot}
\end{equation*}
We denote by $z_a$ be the element of $\mathfrak{z}$ determined by
\begin{equation*}
\langle B(a)u,w \rangle =a([u,w])=\langle J_{z_a} u,w \rangle
\end{equation*}
Since $B(a)$ is skew symmetric and non-degenerate, the dimension of $\mathfrak{z}^{\bot}$ is even, i.e. $dim\mathfrak{z}^{\bot}=2n$.

For a given $0 \neq a \in \mathfrak{z}^*$, we can choose an orthonormal basis
\begin{equation*}
\{E_1(a),E_2(a),\cdots,E_n(a),\overline{E}_1(a),\overline{E}_2(a),\cdots,\overline{E}_n(a)\}
\end{equation*}
of $\mathfrak{z}^{\bot}$ such that
\begin{equation*}
B(a)E_i(a)=|z_a|J_{\frac{z_a}{|z_a|}}E_i(a)=|a|\overline{E}_i(a)
\end{equation*}
and
\begin{equation*}
B(a)\overline{E}_i(a)=-|a|E_i(a)
\end{equation*}
We set $m=dim \mathfrak{z}$. Throughout this paper we assume that $m>1$. We can choose an orthonormal basis $\{\epsilon_1,\epsilon_2,\cdots,\epsilon_m \}$ of $\mathfrak{z}$ such that $a(\epsilon_1)=|a|,a(\epsilon_j)=0,
j=2,3,\cdots,m$. Then we can denote the element of $\mathfrak{g}$ by
\begin{equation*}
(z,t)=(x,y,t)=\underset{i=1}{\overset{n}{\sum}}(x_i E_i+y_i \overline{E}_i )+\underset{j=1}{\overset{m}{\sum}}t_j \epsilon_j
\end{equation*}
We identify G with its Lie algebra $\mathfrak{g}$ by exponential map. The group law on H-type group $G$ has the form
\begin{equation}
(z,t)(z',t')=(z+z',t+t'+\frac{1}{2}[z,z'])  \label{equ:Law}
\end{equation}
where $[z,z']_j=\langle z,U^jz' \rangle$ for a suitable skew symmetric matrix $U^j,j=1,2,\cdots,m$.

\begin{theorem} \indent G is an H-type group with underlying manifold $\mathbb{R}^{2n+m}$, with the group law  $\eqref{equ:Law}$ and the matrix $U^j,j=1,2,\cdots, m$ satisfies the following conditions:\\
$(i)$ $U^j$ is a $2n \times 2n$ skew symmetric and orthogonal matrix, $j=1,2,\cdots, m$.\\
$(ii)$ $U^i U^j+U^j U^i=0,i,j=1,2,\cdots,m$ with $i \neq j$.
\end{theorem}
{\bf Proof.} See \cite{BU}.

\begin{remark}In particular, $\langle z,U^1 z'\rangle=\underset{j=1}{\overset{n}{\sum}}(x_j'y_j-y_j'x_j)$.
\end{remark}

\begin{remark} All the above expressions depend on a given $0 \neq a \in \mathfrak{z}^*$, but we will suppress $a$ from them for simplification.
\end{remark}

\begin{remark}\label{remark}It is well know that H-type algebras are closely related to Clifford modules \cite{R}. H-type algebras can be classified by the standard theory of Clifford algebras. Specially, on H-type group $G$, there is a relation between the dimension of the center and its orthogonal complement space. That is $m+1\leq 2n$ (see \cite{KR}).
\end{remark}

The left invariant vector fields which agree respectively with $\frac{\partial}{\partial x_j},\frac{\partial}{\partial y_j}$ at the origin are given by
\begin{align*}
X_j&=\frac{\partial}{\partial x_j}+\frac{1}{2}\underset{k=1}{\overset{m}{\sum}} \left( \underset{l=1}{\overset{2n}{\sum}}z_l U_{l,j}^k \right) \frac{\partial}{\partial t_k}\\
Y_j&=\frac{\partial}{\partial y_j}+\frac{1}{2}\underset{k=1}{\overset{m}{\sum}} \left( \underset{l=1}{\overset{2n}{\sum}}z_l U_{l,j+n}^k \right) \frac{\partial}{\partial t_k}\\
\end{align*}
where $z_l=x_l,z_{l+n}=y_l,l=1,2,\cdots,n.$

The vector fields $T_k=\frac{\partial}{\partial t_k},k=1,2,\cdots,m$ correspond to the center of $G$. In terms of these vector fields we introduce the sublaplacian $L$ and full Laplacian $\Delta$ respectively
\begin{align}
L&=-\underset{j=1}{\overset{n}{\sum}}(X_j^2 +Y_j^2)=-\Delta_z +\frac{1}{4} |z|^2 T -\underset{k=1}{\overset{m}{\sum}}\langle z,U^k \nabla_z \rangle T_k  \label{eq:Lap}\\
\Delta&=L+T  \nonumber \\
\end{align}
where
\begin{equation*}
\Delta_z=\underset{j=1}{\overset{2n}{\sum}}\frac{\partial^2}{\partial z_j^2}, T=-\underset{k=1}{\overset{m}{\sum}}\frac{\partial^2}{\partial t_k^2},\nabla_z=(\frac{\partial}{\partial z_1},\frac{\partial}{\partial z_2},\cdots,\frac{\partial}{\partial z_{2n}})^t.\\
\end{equation*}

\section[Restriction operator]{Restriction operator}
\quad First we recall some results about the scaled special Hermite expansion. We refer the reader to \cite{T2} and \cite{T3} for details. Let $\lambda >0$. The twisted Laplacian (or the scaled special Hermite expansion) $L_\lambda$ is defined by
\begin{equation*}
L_\lambda=-\Delta_z +\frac{\lambda^2|z|^2}{4}-i\lambda \underset{j=1}{\overset{n}{\sum}}\left( x_j \frac{\partial}{\partial y_j}-y_j \frac{\partial}{\partial x_j} \right).
\end{equation*}
where we identify $z=x+iy \in \mathbb{C}^n$ with $z=(x,y)\in \mathbb{R}^{2n}$.

For $f,g \in L^1(\mathbb{C}^n)$,we define the $\lambda$-twisted convolution by
\begin{equation*}
f \times_\lambda g=\int_{\mathbb{C}^n} f(z-w) g(w) e^{\frac{i \lambda}{2}Im \,z \cdot \overline{w}} \, dw
\end{equation*}

Set Laguerre function $\varphi_k^\lambda(z)=L_k^{n-1}(\frac{1}{2}\lambda|z|^2)e^{-\frac{1}{4}\lambda|z|^2}$, $k=0,1,2,\cdots$,where $L_k^{n-1}$ is the Laguerre polynomial of type $(n-1)$ and degree k.

For $f\in \mathscr{S}(\mathbb{C}^n)$,we have the scaled special Hermite expansion
\begin{equation}
f(z)=\left (\frac{\lambda}{2\pi} \right)^n \underset{k=0}{\overset{+\infty}{\sum}}f \times_\lambda \varphi_k^\lambda(z) \label{eq:f}
\end{equation}
which is an orthogonal form. We also have
\begin{equation}
\label{eq:Plancherel}
||f||^2=\left (\frac{\lambda}{2\pi} \right)^n\underset{k=0}{\overset{+\infty}{\sum}}||f\times_\lambda \varphi_k^\lambda||^2
\end{equation}

Moreover, $f \times_\lambda \varphi_k^\lambda$ ia an eigenfunction of $L_\lambda$ with the eigenvalue $(2k+n)\lambda$ and
\begin{equation}
||f \times_\lambda \varphi_k^\lambda||_2 \leq (2k+n)^{n(\frac{1}{p}-\frac{1}{2})-\frac{1}{2}} \lambda^{n(\frac{1}{p}-\frac{3}{2})} ||f||_p ,\ for \, 1\leq p < \frac{6n+2}{3n+4} \label{eq:eatimate}.
\end{equation}
(see \cite{T1})

Now we turn to the expression for the restriction operator. We may identify $\mathfrak{z}^*$ with $\mathfrak{z}$. Therefore, we will write $\langle a,t \rangle$ instead of $a(t)$ for $a\in \mathfrak{z}^*$ and $t \in \mathfrak{z}$.

\begin{lemma}\label{good}\indent Let $0\neq a \in \mathfrak{z}^*$. If $f(z,t)=e^{-i\langle a,t \rangle}\varphi(z)$, then
\begin{equation*}
Lf(z,t)=e^{-i\langle a,t \rangle}L_{|a|}\varphi(z).
\end{equation*}
\end{lemma}

{\bf Proof.}\indent Because of $\langle a,t \rangle=|a|t_1$ and $\langle z,U^1 \nabla_z \rangle =\underset{j=1}{\overset{n}{\sum}} \left( y_j \frac{\partial}{\partial x_j}-x_j \frac{\partial}{\partial y_j}\right)$, Lemma 3.1 is easily deduced from the expression \eqref{eq:Lap}.

Set $e^a_k(z,t)=e^{-i\langle a,t \rangle}\varphi_k^{|a|}(z)$. For $f\in \mathscr{S}(G)$, let
\begin{equation*}
f^a(z)=\int_{\mathbb{R}^m}f(z,t)e^{i\langle a,t \rangle}\,dt
\end{equation*}
be the Fourier transform of $f$ with respect to the central variable $t$. It is easy to obtain
\begin{equation}
f*e^a_k(z,t)=e^{-i\langle a,t \rangle}f^a \times_{|a|} \varphi_k^{|a|}(z) \label{eq:twisted}
\end{equation}
Note that $f*e^a_k$ is an eigenfunction of $T$ with the eigenvalue $|a|^2$. Furthermore, it follows from Lemma \ref{good} that $f*e^a_k$ is an eigenfunction of $L$ with the eigenvalue $(2k+n)|a|$. Thus $f*e^a_k$ is a joint eigenfunction of the operator $L$ and $T$.

For a Schwartz function $f$ on H-type group,using the inversion formula for the Fourier transform together with \eqref{eq:f} and \eqref{eq:twisted},we have
\begin{align*}
f(z,t)&=\frac{1}{(2\pi)^m}\int_{\mathbb{R}^m} f^a(z)e^{-i\langle a,t \rangle}\, da\\
      &=\frac{1}{(2\pi)^m}\int_{\mathbb{R}^m} \left(\frac{|a|^n}{(2\pi)^n} \underset{k=0}{\overset{+\infty}{\sum}}f^a \times_{|a|} \varphi_k^{|a|}(z)\right)e^{-i\langle a,t \rangle}\, da\\
      &=\frac{1}{(2\pi)^{n+m}}\int_{\mathbb{R}^m}\overset{\infty}{\underset{k=0}{\sum}}f*e^a_k(z,t)|a|^n \, da \\
      &=\int_0^{+\infty}\biggl( \frac{1}{(2\pi)^{n+m}}\overset{\infty}{\underset{k=0}{\sum}}\lambda^{n+m-1}\int_{S^{m-1}}f*e^{\lambda \tilde{a}}_k(z,t)\,d\sigma(\tilde{a})\biggr)\,d\lambda
\end{align*}

The operators $L$ and $T$ extend to a pair of strongly commuting self-adjoint operators. Therefore, they admit a joint spectral decomposition. By the spectral theorem, we can define the joint functional calculus of L and T. Indeed, given a bounded function $h:\mathbb{R}_+ \times \mathbb{R}_+ \rightarrow \mathbb{R}$, we define
\begin{equation*}
h(L,T)f(z,t)=\overset{+\infty}{\underset{0}{\int}}\biggl(\frac{1}{(2\pi)^{n+m}}\overset{\infty}{\underset{k=0}{\sum}}h((2k+n)\lambda,\lambda^2)\lambda^{n+m-1}\int_{S^{m-1}}f*e^{\lambda \tilde{a}}_k(z,t) \, d \sigma(\tilde{a})\biggr)\, d\lambda
\end{equation*}
Then by \eqref{eq:Plancherel} and a simple calculation, we have
\begin{equation*}
\underset{\mathbb{R}^{2n}}{\int}\underset{\mathbb{R}^{m}}{\int}|h(L,T)f(z,t)|^2 \, dzdt\leq  ||h||^2_{L^{\infty}(\mathbb{R}_+ \times \mathbb{R}_+)} ||f||^2_{L^2(G)}
\end{equation*}
We define $\delta(D)f=\underset{\varepsilon \rightarrow 0^+}{\lim}\frac{1}{2\varepsilon}\chi_{(a-\varepsilon,a+\varepsilon)}(D)f$, and $\{{\delta_\mu (D)}\}_{\mu \in \mathbb{R}_+}$ turns out to be the spectral resolution of $D$. ($D=L$ or $\Delta$)

More generally, with the same techniques we can also define operators of the form $\delta_\mu (h(L,T))$ for a suitable function $h$. We assume $h((2k+n)\lambda,\lambda^2)$ is a strictly monotonic differentiable positive function of $\lambda$ on $\mathbb{R}_+$, with the domain $(A,B)$ where $0 \leq A <B \leq +\infty$. Then for each $\mu \in (A,B)$, the equation $h((2k+n)\lambda,\lambda^2)=\mu$ may be solved for each $k$. We denote the solution by $\lambda=\lambda_k(\mu)$ and denotes $\lambda'_k$ the derivative of $\lambda_k$. Replacing in the integral $\lambda$ with $\mu$, we obtain
\begin{equation*}
h(L,T)f(z,t)=\int_A^B \mu \biggl(\frac{1}{(2\pi)^{n+m}}\overset{\infty}{\underset{k=0}{\sum}}
\lambda_k^{n+m-1}(\mu) |\lambda_k'(\mu)| \int_{S^{m-1}}f*e^{\lambda_k (\mu) \tilde{a}}_k(z,t) \, d \sigma(\tilde{a})\biggr)\, d\mu
\end{equation*}
which is the spectral decomposition of $h(L,T)$.

Thus, given a Schwartz function $f$, the spectral decomposition with respect to $h(L,T)$ is
\begin{equation*}
f(z,t)=\int_A^B \biggl(\frac{1}{(2\pi)^{n+m}}\overset{\infty}{\underset{k=0}{\sum}}
\lambda_k^{n+m-1}(\mu) |\lambda_k'(\mu)| \int_{S^{m-1}}f*e^{\lambda_k (\mu) \tilde{a}}_k(z,t) \, d \sigma(\tilde{a})\biggr)\, d\mu
\end{equation*}
Then the spectral resolution of a Schwartz function $f$ is given in terms of the distributions
\begin{equation*}
\mathcal{P}_\mu ^h f(z,t)=\delta_\mu (h(L,T))f(z,t)=\frac{1}{(2\pi)^{n+m}}\overset{\infty}{\underset{k=0}{\sum}}
\lambda_k^{n+m-1}(\mu) |\lambda_k'(\mu)| \int_{S^{m-1}}f*e^{\lambda_k (\mu) \tilde{a}}_k(z,t) \, d \sigma(\tilde{a})
\end{equation*}
Specially, for the full Laplacian $\Delta$, $h(\xi,\eta)=\xi + \eta$, thus we have $\mu=(2k+n)\lambda+\lambda^2$, which yields
 \begin{equation*}
 \lambda_k (\mu)=\frac{1}{2} \sqrt{4\mu +(2k+n)^2} -\frac{2k+n}{2} \quad \text{and} \quad \lambda'_k(\mu)=\frac{1}{\sqrt{4\mu +(2k+n)^2}}
 \end{equation*}
Therefore,
\begin{equation*}
\mathcal{P}_\mu ^\Delta f(z,t)=\frac{1}{(2\pi)^{n+m}}\overset{\infty}{\underset{k=0}{\sum}}
\lambda_k^{n+m-1}(\mu) \lambda_k'(\mu) \int_{S^{m-1}}f*e^{\lambda_k(\mu) \tilde{a}}_k(z,t) \, d \sigma(\tilde{a})\\
\end{equation*}

\section[Restriction Theorem]{Restriction Theorem}
\quad Our main result is the following theorem.
\begin{theorem}\label{main}\indent Let $G$ be an H-type group with the underlying manifold $\mathbb{R}^{2n+m}$, where $m>1$ is the dimension of the center. Let $h(\xi,\eta)=\xi^\alpha+\eta^\beta$, $\alpha,\beta>0$. Then for $1\leq p \leq \frac{2m+2}{m+3}$, we have
if $\alpha<2\beta$,
\begin{equation*}
||\mathcal{P}_\mu ^h f||_{p'} \leq C\mu^{\frac{2}{\alpha}(n+\frac{\alpha}{2\beta}m)(\frac{1}{p}-\frac{1}{2})-1} ||f||_p,\, \mu>1
\end{equation*}
and
\begin{equation*}
||\mathcal{P}_\mu ^h f||_{p'} \leq C\mu^{\frac{2}{\alpha}(n+m)(\frac{1}{p}-\frac{1}{2})-1} ||f||_p,\, 0<\mu\leq1
\end{equation*}
if $\alpha>2\beta$,
\begin{equation*}
||\mathcal{P}_\mu ^h f||_{p'} \leq C\mu^{\frac{2}{\alpha}(n+m)(\frac{1}{p}-\frac{1}{2})-1} ||f||_p,\, \mu>1
\end{equation*}
and
\begin{equation*}
||\mathcal{P}_\mu ^h f||_{p'} \leq C\mu^{\frac{2}{\alpha}(n+\frac{\alpha}{2\beta}m)(\frac{1}{p}-\frac{1}{2})-1} ||f||_p,\, 0<\mu\leq1
\end{equation*}
if $\alpha=2\beta$,
\begin{equation*}
||\mathcal{P}_\mu ^h f||_{p'} \leq C\mu^{\frac{2}{\alpha}(n+m)(\frac{1}{p}-\frac{1}{2})-1} ||f||_p,\, 0<\mu<+\infty
\end{equation*}
\end{theorem}
hold for all Schwartz functions $f$.

First we prove the following abstract statement.
\begin{proposition} \label{proposition}\indent $h((2k+n)\lambda,\lambda^2)$ is a strictly monotonic differentiable positive function of $\lambda$ on $\mathbb{R}_+$, with the domain $(A,B)$ where $0 \leq A <B \leq +\infty$.\\
Then for $1\leq p \leq \frac{2m+2}{m+3}$, the estimate holds
\begin{equation*}
||\mathcal{P}_\mu ^h f||_{p'} \leq C_\mu ||f||_p
\end{equation*}
where
\begin{equation}
\label{equ:Estimate}
C_\mu \leq C \overset{\infty}{\underset{k=0}{\sum}}(2k+n)^{2n(\frac{1}{p}-\frac{1}{2})-1} \lambda_k^{2(n+m)(\frac{1}{p}-\frac{1}{2})-1}(\mu) |\lambda_k'(\mu)|
\end{equation}
for all Schwartz functions $f$ and all positive $\mu \in (A,B)$.
\end{proposition}

{\bf Proof.} \indent Because the scaled special Hermite expansion is orthogonal, we have
\begin{align*}
\left \langle \mathcal{P}_\mu^h f,g \right \rangle
&=\underset{G}{\int}\mathcal{P}_\mu^h f(z,t) \overline{g(z,t)}\, dzdt   \\
&=\underset{G}{\int}\biggl(\frac{1}{(2\pi)^{n+m}}\overset{\infty}{\underset{k=0}{\sum}} \lambda_k^{n+m-1}(\mu) |\lambda_k'(\mu)| \int_{S^{m-1}}f*e^{\lambda_k(\mu) \tilde{a}}_k(z,t) \, d \sigma(\tilde{a})\biggr)
                               \overline{g(z.t)} \,dzdt  \\
&=\frac{1}{(2\pi)^{n+m}}\overset{\infty}{\underset{k=0}{\sum}}\lambda_k^{n+m-1}(\mu) |\lambda_k'(\mu)|
                              \int_{S^{m-1}} \underset{\mathbb{R}^{2n}}{\int} f^{\lambda_k(\mu) \tilde{a}} \times_{\lambda_k(\mu)} \varphi_k ^{\lambda_k(\mu)} (z) \overline{g^{\lambda_k(\mu) \tilde{a}} (z)} \,dz d\sigma(\tilde{a}) \\
&=\frac{1}{(2\pi)^{2n+m}}\overset{\infty}{\underset{k=0}{\sum}}\lambda_k^{2n+m-1}(\mu) |\lambda_k'(\mu)|
                              \int_{S^{m-1}} \underset{\mathbb{R}^{2n}}{\int} f^{\lambda_k(\mu) \tilde{a}} \times_{\lambda_k(\mu)} \varphi_k ^{\lambda_k(\mu)} (z) \\
&\quad \quad \quad \quad \quad \quad \quad \quad \quad  \times\overline{g^{\lambda_k(\mu) \tilde{a}} \times_{\lambda_k(\mu)} \varphi_k ^{\lambda_k(\mu)} (z)} \,dz d\sigma(\tilde{a}) \\
 & \leq \frac{1}{(2\pi)^{2n+m}}\overset{\infty}{\underset{k=0}{\sum}}\lambda_k^{2n+m-1}(\mu) |\lambda_k'(\mu)|
                              \int_{S^{m-1}}||f^{\lambda_k(\mu) \tilde{a}} \times_{\lambda_k(\mu)} \varphi_k ^{\lambda_k(\mu)}||_2 \\
 &\quad \quad \quad \quad \quad \quad \quad \quad \quad \times
                              ||g^{\lambda_k(\mu) \tilde{a}} \times_{\lambda_k(\mu)} \varphi_k ^{\lambda_k(\mu)}||_2 \, d\sigma(\tilde{a}) \\
\end{align*}
Because of \eqref{eq:twisted}, we have
\begin{equation}
\label{equ:1.1}
||f^{\lambda_k(\mu) \tilde{a}} \times_{\lambda_k(\mu)} \varphi_k ^{\lambda_k(\mu)}||_2 \leq C(2k+n)^{n(\frac{1}{p}-\frac{1}{2})-\frac{1}{2}} \lambda_k^{n(\frac{1}{p}-\frac{3}{2})}(\mu)||f^{\lambda_k(\mu) \tilde{a}}||_p,1\leq p < \frac{6n+2}{3n+4}
\end{equation}
From Remark \ref{remark}, we have $m+1 \leq 2n$. Thus $\frac{2m+2}{m+3} <\frac{6n+2}{3n+4}$.
Therefore, by $(\ref{equ:1.1})$, applying H$\ddot{o}$lder inequality and the Minkowski inequality, we get
\begin{align*}
&\quad \left \langle \mathcal{P}_\mu^h f,g \right \rangle\\
&\leq C\overset{\infty}{\underset{k=0}{\sum}}(2k+n)^{2n(\frac{1}{p}-\frac{1}{2})-1}
\lambda_k^{2n(\frac{1}{p}-\frac{1}{2})+m-1}(\mu) |\lambda_k'(\mu)| \int_{S^{m-1}}||f^{\lambda_k(\mu) \tilde{a}}||_p ||g^{\lambda_k(\mu) \tilde{a}}||_p\,d\sigma(\tilde{a})\\
&\leq C\overset{\infty}{\underset{k=0}{\sum}}(2k+n)^{2n(\frac{1}{p}-\frac{1}{2})-1} \lambda_k^{2n(\frac{1}{p}-\frac{1}{2})+m-1}(\mu) |\lambda_k'(\mu)| \left( \underset{S^{m-1}}{\int} \Big(\underset{\mathbb{R}^{2n}}{\int}|f^{\lambda_k(\mu) \tilde{a}}(z)|^p \, dz \Big) ^{\frac{2}{p}}d\sigma(\tilde{a}) \right) ^{\frac{1}{2}} \\
& \quad \quad \quad \quad \quad \quad \quad \quad \times \left( \underset{S^{m-1}}{\int} \Big(\underset{\mathbb{R}^{2n}}{\int}|g^{\lambda_k(\mu) \tilde{a}}(z)|^p \, dz \Big) ^{\frac{2}{p}}d\sigma(\tilde{a}) \right) ^{\frac{1}{2}}\\
&\leq C\overset{\infty}{\underset{k=0}{\sum}}(2k+n)^{2n(\frac{1}{p}-\frac{1}{2})-1} \lambda_k^{2n(\frac{1}{p}-\frac{1}{2})+m-1}(\mu) |\lambda_k'(\mu)| \biggl( \underset{\mathbb{R}^{2n}}{\int} \bigl(\underset{S^{m-1}}{\int}|f^{\lambda_k(\mu) \tilde{a}}(z)|^2 \, d\sigma(\tilde{a}) \bigr) ^{\frac{p}{2}} dz \biggr) ^{\frac{1}{p}} \\
& \quad \quad \quad \quad \quad \quad \quad \quad \times \left( \underset{\mathbb{R}^{2n}}{\int} \Big(\underset{S^{m-1}}{\int}|g^{\lambda_k(\mu) \tilde{a}}(z)|^2 \, d\sigma(\tilde{a}) \Big) ^{\frac{p}{2}} dz \right) ^{\frac{1}{p}} \\
\end{align*}
Denote by $f_{k,\Delta}(z,t)=f(z,\frac{t}{\lambda_k(\mu)})$, and then we have
\begin{equation*}
\Big( \underset{S^{m-1}}{\int}|f^{\lambda_k(\mu) \tilde{a}}(z)|^2 \, d\sigma(\tilde{a}) \Big) ^{\frac{1}{2}}=\lambda_k(\mu)^{-m}\Big( \underset{S^{m-1}}{\int}|f_{k,\Delta}^{\tilde{a}}(z)|^2 \, d\sigma(\tilde{a}) \Big) ^{\frac{1}{2}}
\end{equation*}
It follows from Stein-Tomas restriction theorem that
\begin{equation*}
\Big( \underset{S^{m-1}}{\int}|f_{k,\Delta}^{\tilde{a}}(z)|^2 \, d\sigma(\tilde{a}) \Big) ^{\frac{1}{2}}
\leq C \Big( \underset{\mathbb{R}^m} {\int}|f_{k,\Delta}(z,t)|^p \, dz \Big) ^{\frac{1}{p}}
=C\lambda_k(\mu)^{\frac{m}{p}}\Big( \underset{\mathbb{R}^m} {\int}|f(z,t)|^p \, dz \Big) ^{\frac{1}{p}}.
\end{equation*}
Therefore,
\begin{equation*}
\left \langle \mathcal{P}_\mu^h f,g \right \rangle \leq C \overset{\infty}{\underset{k=0}{\sum}}(2k+n)^{2n(\frac{1}{p}-\frac{1}{2})-1} \lambda_k^{2(n+m)(\frac{1}{p}-\frac{1}{2})-1}(\mu) |\lambda_k'(\mu)| ||f||_p ||g||_p
\end{equation*}
proves the statement.

To obtain Theorem \ref{main}, it suffices to show the convergence of the series in \eqref{equ:Estimate}. Next we will exploit the following estimates, which can be easily proved by comparing the sums with the corresponding integrals:
\begin{lemma}\label{sum}
Fix $\nu \in \mathbb{R}$. There exists $C_\nu>0$ such that for $A>0$ and $n\in\mathbb{Z}_+$, we have
\begin{align}
\underset{2m+n\geq A}{\underset{m\in\mathbb{N}}{\sum}}(2m+n)^\nu &\leq C_\nu A^{\nu+1}, \quad \nu<-1;\label{sum1}\\
\underset{2m+n\leq A}{\underset{m\in\mathbb{N}}{\sum}}(2m+n)^\nu &\leq C_\nu A^{\nu+1}, \quad \nu>-1.\label{sum2}
\end{align}
\end{lemma}

Now theorem \ref{main} follows from the result in the following lemma.
\begin{lemma} \label{lemma}\indent Let $h(\xi,\eta)=\xi^\alpha+\eta^\beta$, $\alpha,\beta>0$. The series in \eqref{equ:Estimate} has the estimate\\
if $\alpha<2\beta$,
\begin{equation*}
C_\mu \leq C\mu^{\frac{2}{\alpha}(n+\frac{\alpha}{2\beta}m)(\frac{1}{p}-\frac{1}{2})-1},\, \mu>1
\end{equation*}
and
\begin{equation*}
C_\mu\leq C\mu^{\frac{2}{\alpha}(n+m)(\frac{1}{p}-\frac{1}{2})-1},\, 0<\mu\leq1
\end{equation*}
if $\alpha>2\beta$,
\begin{equation*}
C_\mu \leq C\mu^{\frac{2}{\alpha}(n+m)(\frac{1}{p}-\frac{1}{2})-1},\, \mu>1
\end{equation*}
and
\begin{equation*}
C_\mu \leq C\mu^{\frac{2}{\alpha}(n+\frac{\alpha}{2\beta}m)(\frac{1}{p}-\frac{1}{2})-1},\, 0<\mu\leq1
\end{equation*}
if $\alpha=2\beta$,
\begin{equation*}
C_\mu \leq C\mu^{\frac{2}{\alpha}(n+m)(\frac{1}{p}-\frac{1}{2})-1},\, 0<\mu<+\infty
\end{equation*}
\end{lemma}

{\bf Proof.}\indent $h(\xi,\eta)=\xi^\alpha+\eta^\beta$, $\alpha,\beta>0$, thus we have $\mu=(2k+n)^\alpha \lambda_k^\alpha(\mu)+\lambda_k^{2\beta}(\mu)$, which yields
 \begin{equation*}
 \lambda_k'(\mu)=\frac{1}{\alpha (2k+n)^\alpha \lambda_k^{\alpha-1}(\mu)+2\beta\lambda_k^{2\beta-1}(\mu)}
 \end{equation*}
 To study the convergence of this series, we need to distinguish three cases according to the relative of $\alpha$ and $2\beta$:$\alpha<2\beta$,$\alpha>2\beta$ and $\alpha=2\beta$. In order not to burden the exposition, we only prove the case $\alpha<2\beta$, and the other cases are analogous.\\
 If $\alpha<2\beta$:\\
 when $\mu \leq 1$, it is easy to see that $\lambda_k(\mu)\sim \frac{\mu^{\frac{1}{\alpha}}}{2k+n}$ and $\lambda_k'(\mu)\sim \frac{\mu^{\frac{1}{\alpha}-1}}{2k+n}$, so that the series
\begin{equation}
\label{equ:A}
\begin{split}
C_\mu
      &\leq C \overset{\infty}{\underset{k=0}{\sum}}(2k+n)^{2n(\frac{1}{p}-\frac{1}{2})-1} \lambda_k^{2(n+m)(\frac{1}{p}-\frac{1}{2})-1}(\mu) |\lambda_k'(\mu)|\\
      &\leq C\overset{\infty}{\underset{k=0}{\sum}}(2k+n)^{2n(\frac{1}{p}-\frac{1}{2})-1} \biggl(\frac{\mu^{\frac{1}{\alpha}}}{2k+n} \biggr)^{2(n+m)(\frac{1}{p}-\frac{1}{2})-1} \frac{\mu^{\frac{1}{\alpha}-1}}{2k+n}\\
      &\leq C\mu^{\frac{2}{\alpha}(n+m)(\frac{1}{p}-\frac{1}{2})-1}\overset{\infty}{\underset{k=0}{\sum}}\frac{1}{(2k+n)^{2m(\frac{1}{p}-\frac{1}{2})+1}}\\
      &\leq C\mu^{\frac{2}{\alpha}(n+m)(\frac{1}{p}-\frac{1}{2})-1} \\
\end{split}
\end{equation}
converges.

When $\mu>1$, we split the sum into two parts, the sum over those $k$ such that $(2k+n)^\alpha \lambda_k^\alpha(\mu)\geq \lambda^{2\beta}(\mu)$ and those such that $(2k+n)^\alpha \lambda_k^\alpha(\mu)<\lambda^{2\beta}(\mu)$. They are denoted by $\uppercase \expandafter {\romannumeral 1}$ and $\uppercase \expandafter {\romannumeral 2}$ respectively.

For the first part, $(2k+n)^\alpha \lambda_k^\alpha(\mu)\geq \lambda^{2\beta}(\mu)$ implies
\begin{equation*}
\lambda_k(\mu)\sim \frac{\mu^{\frac{1}{\alpha}}}{2k+n},\lambda_k'(\mu)\sim \frac{\mu^{\frac{1}{\alpha}-1}}{2k+n},\text{and} \quad 2k+n\geq \mu^{\frac{2\beta-\alpha}{2\alpha \beta}}
\end{equation*}
Then we control the first part $\uppercase \expandafter {\romannumeral 1}$ by
\begin{align*}
\uppercase \expandafter {\romannumeral 1}
                         &\leq C \underset{2k+n\geq \mu^{\frac{2\beta-\alpha}{2\alpha \beta}}}{\sum} (2k+n)^{2n(\frac{1}{p}-\frac{1}{2})-1} \lambda_k^{2(n+m)(\frac{1}{p}-\frac{1}{2})-1}(\mu) |\lambda_k'(\mu)|\\
                         &\leq C \underset{2k+n\geq \mu^{\frac{2\beta-\alpha}{2\alpha \beta}}}{\sum}(2k+n)^{2n(\frac{1}{p}-\frac{1}{2})-1} \biggl(\frac{\mu^{\frac{1}{\alpha}}}{2k+n} \biggr)^{2(n+m)(\frac{1}{p}-\frac{1}{2})-1}
                         \frac{\mu^{\frac{1}{\alpha}-1}}{2k+n}\\
                         &\leq C \mu^{\frac{2}{\alpha}(n+m)(\frac{1}{p}-\frac{1}{2})-1}\underset{2k+n\geq \mu^{\frac{2\beta-\alpha}{2\alpha \beta}}}{\sum}\frac{1}{(2k+n)^{2m(\frac{1}{p}-\frac{1}{2})+1}}
\end{align*}
Then by \eqref{sum1}, we have
\begin{equation}
\uppercase \expandafter {\romannumeral 1}\leq C \mu^{\frac{2}{\alpha}(n+m)(\frac{1}{p}-\frac{1}{2})-1} \frac{1}{(\mu^{\frac{2\beta-\alpha}{2\alpha \beta}})^{2m(\frac{1}{p}-\frac{1}{2})}}\leq C\mu^{\frac{2}{\alpha}(n+\frac{\alpha}{2\beta}m)(\frac{1}{p}-\frac{1}{2})-1} \label{equ:B}
\end{equation}
For the second part, $(2k+n)^\alpha \lambda_k^\alpha(\mu)<\lambda^{2\beta}(\mu)$ implies
\begin{equation*}
\lambda_k(\mu)\sim \mu^{\frac{1}{2\beta}},\lambda_k'(\mu)\sim \mu^{\frac{1}{2\beta}-1},\text{and} \quad 2k+n < \mu^{\frac{2\beta-\alpha}{2\alpha \beta}}
\end{equation*}
Then we control the second part $\uppercase \expandafter {\romannumeral 2}$ by
\begin{align*}
\uppercase \expandafter {\romannumeral 2}
                         &\leq C \underset{2k+n<\mu^{\frac{2\beta-\alpha}{2\alpha \beta}}}{\sum} (2k+n)^{2n(\frac{1}{p}-\frac{1}{2})-1} \lambda_k^{2(n+m)(\frac{1}{p}-\frac{1}{2})-1}(\mu) |\lambda_k'(\mu)|\\
                         &\leq C \underset{2k+n<\mu^{\frac{2\beta-\alpha}{2\alpha \beta}}}{\sum}(2k+n)^{2n(\frac{1}{p}-\frac{1}{2})-1} \biggl(\mu^{\frac{1}{2\beta}}\biggr)^{2(n+m)(\frac{1}{p}-\frac{1}{2})-1}
                         \mu^{\frac{1}{2\beta}-1}\\
                         &\leq C \mu^{\frac{1}{\beta}(n+m)(\frac{1}{p}-\frac{1}{2})-1}\underset{2k+n<\mu^{\frac{2\beta-\alpha}{2\alpha \beta}}}{\sum}(2k+n)^{2n(\frac{1}{p}-\frac{1}{2})-1} \\
\end{align*}
Because $1\leq p \leq \frac{2m+2}{m+3}$ , we obtain $2n(\frac{1}{p}-\frac{1}{2})-1 \geq -1$. Hence, by \eqref{sum2} we get
\begin{equation*}
\underset{2k+n<\mu^{\frac{2\beta-\alpha}{2\alpha \beta}}}{\sum}(2k+n)^{2n(\frac{1}{p}-\frac{1}{2})-1}
\lesssim (\mu^{\frac{2\beta-\alpha}{2\alpha \beta}})^{2n(\frac{1}{p}-\frac{1}{2})}=\mu^{(\frac{2}{\alpha}-\frac{1}{\beta})n(\frac{1}{p}-\frac{1}{2})}
\end{equation*}
Thus, for the second part we also have
\begin{equation}
\uppercase \expandafter {\romannumeral 2}\leq C\mu^{\frac{2}{\alpha}(n+\frac{\alpha}{2\beta}m)(\frac{1}{p}-\frac{1}{2})-1} \label{equ:C}
\end{equation}
Finally, the estimate for case $\alpha < 2\beta$ follows from \eqref{equ:A}, \eqref{equ:B} and \eqref{equ:C}. This completes the proof of the first case.

Combining Proposition \ref{proposition} and Lemma \ref{lemma}, Theorem \ref{main} comes out easily.

Specially, in case of $\Delta=L+T$, $h(\xi,\eta)=\xi+\eta$, we obtain the restriction theorem associated with the full Laplacian on H-type groups.
\begin{corollary}\indent For $1\leq p \leq \frac{2m+2}{m+3}$, the estimates
\begin{equation*}
||\mathcal{P}_\mu ^\Delta f||_{p'} \leq C \mu^{(2n+m)(\frac{1}{p}-\frac{1}{2})-1} ||f||_p , \,\mu>1
\end{equation*}
and
\begin{equation*}
||\mathcal{P}_\mu ^\Delta f||_{p'} \leq C \mu^{2(n+m)(\frac{1}{p}-\frac{1}{2})-1} ||f||_p ,\, 0<\mu \leq 1
\end{equation*}
hold for all Schwartz functions $f$.
\end{corollary}

Similarly to what we have done so far in Theorem \ref{main}, we now discuss other operators with the form of the joint functional calculus of $L$ and $T$. We obtain the following results. We omit the arguments which are really similar to that of Theorem \ref{main}.

\begin{theorem}\indent Let $h(\xi,\eta)=(\xi^\alpha+\eta^\beta)^{-1}$, $\alpha,\beta>0$. Then for $1\leq p \leq \frac{2m+2}{m+3}$, we have\\
if $\alpha<2\beta$,
\begin{equation*}
||\mathcal{P}_\mu ^h f||_{p'} \leq C\mu^{-\frac{2}{\alpha}(n+m)(\frac{1}{p}-\frac{1}{2})-1} ||f||_p,\, \mu>1
\end{equation*}
and
\begin{equation*}
||\mathcal{P}_\mu ^h f||_{p'} \leq C\mu^{-\frac{2}{\alpha}(n+\frac{\alpha}{2\beta}m)(\frac{1}{p}-\frac{1}{2})-1} ||f||_p,\, 0<\mu\leq1
\end{equation*}
if $\alpha>2\beta$,
\begin{equation*}
||\mathcal{P}_\mu ^h f||_{p'} \leq C\mu^{-\frac{2}{\alpha}(n+\frac{\alpha}{2\beta}m)(\frac{1}{p}-\frac{1}{2})-1} ||f||_p,\, \mu>1
\end{equation*}
and
\begin{equation*}
||\mathcal{P}_\mu ^h f||_{p'} \leq C\mu^{-\frac{2}{\alpha}(n+m)(\frac{1}{p}-\frac{1}{2})-1} ||f||_p,\,0<\mu\leq1
\end{equation*}
if $\alpha=2\beta$,
\begin{equation*}
||\mathcal{P}_\mu ^h f||_{p'} \leq C\mu^{-\frac{2}{\alpha}(n+m)(\frac{1}{p}-\frac{1}{2})-1} ||f||_p,\, 0<\mu<+\infty
\end{equation*}
hold for all Schwartz functions $f$.
\end{theorem}

\begin{theorem}\indent Let $h(\xi, \eta)=(1+\xi)^{-1}$. For $1\leq p \leq \frac{2m+2}{m+3}$, the estimates
\begin{equation*}
||\mathcal{P}_\mu^h f||_{p'} \leq C \mu^{-2(n+m)(\frac{1}{p}-\frac{1}{2})-1} ||f||_p , \,\text{when $\mu \rightarrow 0^+$}
\end{equation*}
and
\begin{equation*}
||\mathcal{P}_\mu ^h f||_{p'} \leq C (1-\mu)^{2(n+m)(\frac{1}{p}-\frac{1}{2})-1} ||f||_p ,\, \text{when $\mu \rightarrow 1^-$}
\end{equation*}
hold for all Schwartz functions $f$.
\end{theorem}

More generally, we have
\begin{theorem}\indent Let $h(\xi,\eta)=(\xi^\alpha+\eta^\beta)^\gamma$, $\alpha, \beta, \gamma>0$. Then for $1\leq p \leq \frac{2m+2}{m+3}$, we have\\
if $\alpha<2\beta$,
\begin{equation*}
||\mathcal{P}_\mu ^h f||_{p'} \leq C\mu^{\frac{2}{\alpha \gamma}(n+\frac{\alpha}{2\beta}m)(\frac{1}{p}-\frac{1}{2})-1} ||f||_p,\,  \mu>1
\end{equation*}
and
\begin{equation*}
||\mathcal{P}_\mu ^h f||_{p'} \leq C\mu^{\frac{2}{\alpha \gamma}(n+m)(\frac{1}{p}-\frac{1}{2})-1} ||f||_p,\, 0<\mu\leq1
\end{equation*}
if $\alpha>2\beta$,
\begin{equation*}
||\mathcal{P}_\mu ^h f||_{p'} \leq C\mu^{\frac{2}{\alpha \gamma}(n+m)(\frac{1}{p}-\frac{1}{2})-1} ||f||_p,\, \mu>1
\end{equation*}
and
\begin{equation*}
||\mathcal{P}_\mu ^h f||_{p'} \leq C\mu^{\frac{2}{\alpha \gamma}(n+\frac{\alpha}{2\beta}m)(\frac{1}{p}-\frac{1}{2})-1} ||f||_p,\,0<\mu\leq1
\end{equation*}
if $\alpha=2\beta$,
\begin{equation*}
||\mathcal{P}_\mu ^h f||_{p'} \leq C\mu^{\frac{2}{\alpha \gamma}(n+m)(\frac{1}{p}-\frac{1}{2})-1} ||f||_p,\, 0<\mu<+\infty
\end{equation*}
\end{theorem}
hold for all Schwartz functions $f$.\\

\begin{theorem}\indent Let $h(\xi,\eta)=(\xi^\alpha+\eta^\beta)^{-\gamma}$, $\alpha,\beta,\gamma>0$. Then for $1\leq p \leq \frac{2m+2}{m+3}$, we have\\
if $\alpha<2\beta$,
\begin{equation*}
||\mathcal{P}_\mu ^h f||_{p'} \leq C\mu^{-\frac{2}{\alpha \gamma}(n+m)(\frac{1}{p}-\frac{1}{2})-1} ||f||_p,\, \mu>1
\end{equation*}
and
\begin{equation*}
||\mathcal{P}_\mu ^h f||_{p'} \leq C\mu^{-\frac{2}{\alpha \gamma}(n+\frac{\alpha}{2\beta}m)(\frac{1}{p}-\frac{1}{2})-1} ||f||_p,\,0<\mu\leq1
\end{equation*}
if $\alpha>2\beta$,
\begin{equation*}
||\mathcal{P}_\mu ^h f||_{p'} \leq C\mu^{-\frac{2}{\alpha \gamma}(n+\frac{\alpha}{2\beta}m)(\frac{1}{p}-\frac{1}{2})-1} ||f||_p,\, \mu>1
\end{equation*}
and
\begin{equation*}
||\mathcal{P}_\mu ^h f||_{p'} \leq C\mu^{-\frac{2}{\alpha \gamma}(n+m)(\frac{1}{p}-\frac{1}{2})-1} ||f||_p,\, 0<\mu\leq1
\end{equation*}
if $\alpha=2\beta$,
\begin{equation*}
||\mathcal{P}_\mu ^h f||_{p'} \leq C\mu^{-\frac{2}{\alpha \gamma}(n+m)(\frac{1}{p}-\frac{1}{2})-1} ||f||_p,\, 0<\mu<+\infty
\end{equation*}
hold for all Schwartz functions $f$.
\end{theorem}

\begin{theorem}\indent Let $h(\xi, \eta)=(1+\xi^\alpha +\eta^\beta)^{-\gamma}$, $\alpha, \beta, \gamma>0$. For $1\leq p \leq \frac{2m+2}{m+3}$, the estimates\\
if $\alpha \leq 2\beta$,
\begin{equation*}
||\mathcal{P}_\mu ^h f||_{p'} \leq C\mu^{-\frac{2}{\alpha \gamma}(n+\frac{\alpha}{2\beta}m)(\frac{1}{p}-\frac{1}{2})-1} ||f||_p,\, \text{when $\mu \rightarrow 0^+$}
\end{equation*}
and
\begin{equation*}
||\mathcal{P}_\mu ^h f||_{p'} \leq C(1-\mu^{\frac{1}{\gamma}})^{\frac{2}{\alpha}(n+m)(\frac{1}{p}-\frac{1}{2})-1} ||f||_p,\, \text{when $\mu \rightarrow 1^-$}
\end{equation*}
if $\alpha>2\beta$,
\begin{equation*}
||\mathcal{P}_\mu ^h f||_{p'} \leq C\mu^{-\frac{2}{\alpha \gamma}(n+m)(\frac{1}{p}-\frac{1}{2})-1} ||f||_p,\, \text{when $\mu \rightarrow 0^+$}
\end{equation*}
and
\begin{equation*}
||\mathcal{P}_\mu ^h f||_{p'} \leq C(1-\mu^{\frac{1}{\gamma}})^{\frac{2}{\alpha}(n+\frac{\alpha}{2\beta}m)(\frac{1}{p}-\frac{1}{2})-1} ||f||_p,\, \text{when $\mu \rightarrow 1^-$}
\end{equation*}
hold for all Schwartz functions $f$.
\end{theorem}

\section[Sharpness of the range p]{Sharpness of the range $p$}
\quad In this section we only give an example to show that the range of $p$ in the restriction theorem associated with the full Laplacian $\Delta$ is sharp. The example is constructed similarly to the counterexample of M$\ddot{u}$ller \cite{M}, which shows that the estimates between Lebesgue spaces for the operators $\mathcal{P}^\Delta_\mu$ are necessarily trivial.

Let $\varphi \in C_c^{\infty}(\mathbb{R}^m)$ be a radial function, such that $\varphi(a)=\psi(|a|)$, where $\psi \in C_c^{\infty}(\mathbb{R})$, $\psi=1$ on a neighborhood of the point $n$ and $\psi=0$ near 0. Let $h$ be a Schwartz function on $\mathbb{R}^m$ and define
\begin{equation*}
f(z,t)=\int_{\mathbb{R}^m} \varphi(a)\hat{h}(a) e^{-\frac{|a|}{4}|z|^2} e^{-i\langle a,t\rangle} |a|^n \, da
\end{equation*}

Denote $g(z,t)=\int_{\mathbb{R}^m} \varphi(a) e^{-\frac{|a|}{4}|z|^2} e^{-i\langle a,t\rangle} |a|^n \, da
           =\int_{\mathbb{R}^{m+2n}}\varphi(a) e^{-\frac{|\xi|^2}{|a|}} e^{-i(\langle a,t \rangle+\langle \xi,z\rangle)} \, d\xi da$.

Hence $\widehat{g(\xi,a)}=\varphi(a) e^{-\frac{|\xi|^2}{|a|}}$, which shows that $\hat{g}$ and consequently $g$ are Schwartz functions. On the other hand, we have $f=h*_tg$, where $"*_t"$ denotes the involution about the central variable. By \text{Lemma} \ref{good}, we have  $\Delta \left(e^{-i\langle a,t\rangle} e^{-\frac{|a|}{4}|z|^2}\right)=(n\lambda+\lambda^2)e^{-i\langle a,t\rangle} e^{-\frac{|a|}{4}|z|^2}$. Therefore, we write $f$ by the integration with polar coordinates as
\begin{align*}
f(z,t)&=\int_0^{+\infty}\biggl(\lambda^{n+m-1}\psi(\lambda) e^{-\frac{\lambda}{4}|z|^2}\underset{S^{m-1}}{\int} \hat{h}(\lambda w)e^{-i \lambda
      \langle w,t\rangle} \, d\sigma(w) \biggr)d\lambda \\
      &=\int_0^{+\infty}\biggl(\lambda_\Delta(\mu)^{n+m-1}\lambda'_\Delta(\mu)\psi\big(\lambda_\Delta(\mu)\big) e^{-\frac{\lambda_\Delta(\mu)}{4}|z|^2}\underset{S^{m-1}}{\int} \hat{h}\big(\lambda_\Delta(\mu) w\big)e^{-i \lambda_\Delta(\mu)\langle w,t\rangle} \, d\sigma(w) \biggr)d\mu \\
      &=\int_0^{+\infty} \mathcal{P}^\Delta_\mu f(z,t) \, d\mu\\
\end{align*}
where
\begin{align*}
\mathcal{P}^\Delta_\mu f(z,t)&=\lambda_\Delta(\mu)^{n+m-1}\lambda'_\Delta(\mu)\psi\big(\lambda_\Delta(\mu)\big) e^{-\frac{\lambda_\Delta(\mu)}{4}|z|^2}\underset{S^{m-1}}{\int} \hat{h}\big(\lambda_\Delta(\mu) w\big)e^{-i \lambda_\Delta(\mu)
      \langle w,t\rangle} \, d\sigma(w)\\
\lambda_\Delta(\mu)&=\frac{\sqrt{n^2+4\mu}-n}{2}
\end{align*}

Therefore, specially let $\mu=2n^2$, we have $\lambda_\Delta(2n^2)=n$, $\lambda'_\Delta(2n^2)=\frac{1}{3n}$ and
\begin{equation*}
\begin{split}
\mathcal{P}^\Delta_{2n^2}f(z,t)&=\frac{1}{3}n^{n+m-2}e^{-\frac{n|z|^2}{4}}\underset{S^{m-1}}{\int}\hat{h}(nw)e^{-in\langle w,t \rangle} \, d\sigma(w)\\
                   &=\frac{1}{3}n^{n-1}e^{-\frac{n|z|^2}{4}}h*\widehat{d\sigma_n}(t)\\
\end{split}
\end{equation*}
\indent From the restriction theorem associated the full Laplacian on H-type groups, we have the estimate $||\mathcal{P}^\Delta_{2n^2}f||_{L^{p'}(G)}\leq C ||f||_{L^p(G)}$.

Because of
\begin{equation}
||\mathcal{P}^\Delta_{2n^2}f||_{L^{p'}(G)}=C||h*\widehat{d\sigma_n}||_{L^{p^{'}}(\mathbb{R}^m)}
\end{equation}
 and
\begin{equation}
||f||_{L^p(G)}\leq ||h||_{L^p(\mathbb{R}^m)}||g||_{L^1_t{L^p_z}}\lesssim||h||_{L^p(\mathbb{R}^m)},
\end{equation}
where the mixed Lebesgue norm is defined by
\begin{equation*}
||g||_{L^1_t{L^p_z}}=\biggl( \int_{\mathbb{R}^2n}\biggl( \int_{\mathbb{R}^m} |f(z,t)| \,dt \biggr)^p \, dz\biggr)^{\frac{1}{p}},
\end{equation*}
we have $||h*\widehat{d\sigma_n}||_{L^{p'}(\mathbb{R}^m)} \leq C||h||_{L^p(\mathbb{R}^m)}$.

From the sharpness of Stein-Tomas theorem which is guaranteed by the Knapp counterexample, it would imply $p\leq \frac{2m+2}{m+3}$. Hence the range of $p$ can not be extended. With the same tricks we can also prove the range of $p$ for the restriction theorem associated with the functional calculus is also sharp.\\

{\bf Acknowledgements:}
The work is performed while the second author studies as a joint Ph.D. student in the mathematics department of Christian-Albrechts-Universit$\ddot{a}$t zu Kiel. She is deeply grateful to Professor Detlef M$\ddot{u}$ller for generous discussions and his continuous encouragement.

\begin{flushleft}
\vspace{0.3cm}\textsc{Heping Liu\\School of Mathematical
Sciences\\Peking University\\Beijing 100871\\People's
Republic of China\\}
\emph{E-mail address}: \text{hpliu@pku.edu.cn}

\vspace{0.3cm}\textsc{Manli Song\\School of Mathematical
Sciences\\Peking University\\Beijing 100871\\People's
Republic of China\\}
\emph{E-mail address}: \text{songmanli@pku.edu.cn}
\end{flushleft}

\end{document}